\documentclass{amsart}
\usepackage{graphicx}
\newtheorem{thm}{Theorem}[section]

\theoremstyle{definition}
\newtheorem{defn}[thm]{Definition}
\theoremstyle{remark}

\numberwithin{equation}{section}

\begin{document}

\title[]{Certain Subclasses of Convex functions with positive and missing Coefficients by using a fixed point }%
\author{Sh. Najafzadeh, M. Eshaghi Gordji and A. Ebadian}%

\thanks{}%
\subjclass[2000]{30C45,30C50} \keywords{Subordination; P-valent
function; Coefficient estimate; Fixed point; Distortion bound and
Radii of starlikeness and Convexity.}
\begin{abstract}
By considering a fixed point in unit disk $\Delta$, a new class of
univalent convex functions is defined. Coefficient inequalities,
integral operator and extreme points of this class are obtained.\\

\end{abstract}
\maketitle
\section{Introduction }
Let $w$ be fixed point in $\Delta =\{z:|z|<1\}$ and
$$A(w)=\{f\in
H(\Delta) : f(w)=f'(w)-1=0\}.$$
 Let $M_{w}=\{f\in A(w):$ f is
Univalent in $\Delta\}$, and $M_{w}$ denoted the subclass of $A(w)$
consist of the functions of the form
\begin{equation}
f(z)=\frac{1}{z-w}+\sum^{+\infty}_{n=k}a_{n}(z-w)^{n} , a_{n}\geq0 ,
z\neq w
\end{equation}
 The function $f(z)$ in
$M_{w}$ is said to be Convex of order $\eta(0<\eta<1)$ if and
only if
$$Re\{1+\frac{(z-w)f''(z)}{f'(z)}\}>\eta \ \ \ ,\ \ (z\in\Delta).$$
The similar definitions for uniformly convex and starlike functions
are introduced by Goodman in [$4$]. For the function $f(z)$ in
$M_{w}$, Ghanim and Darus have defined an operator as follows:
 $$ I^{0}f(z)=f(z)$$
$$I^{1}f(z)=(z-w)f'(z)+\frac{2}{z-w}$$
$$I^{2}f(z)=(z-w)(I^{1}f(z))'+\frac{2}{z-w}$$
and for \ \ $ m=1,2,...$
\begin{equation}I^{m}f(z)=(z-w)(I^{m-1}f(z))'+\frac{2}{z-w}\end{equation}
$$=\frac{1}{z-w}+\sum^{+\infty}_{n=k}n^{m}a_{n}(z-w)^{n}.$$
We note that various other subclasses have been studied rather
extensively by many authors. (see also $[3]$ and $[6]$)
\begin{defn}A function $f(z)$ belonging to the class
$M_{w}$ is in the class $M_{w}(A,B,m)$ if it satisfies the
condition
\begin{equation}
 |\frac{\frac{1}{2}[(z-w)\frac{(I^{m}f(z))''}{(I^{m}f(z))'}+1]+\frac{1}{2}}{\frac{B}{2}[(z-w)\frac{(I^{m}f(z))''}{(I^{m}f(z))'}+1]+\frac{A+1}{2}}|\leq1
 \end{equation}
for some $-1\leq B<A<1$ , $0\leq A\leq 1$.
\end{defn}
\section{Main Results }
First we obtain coefficient inequalities for functions in
$M_{w}(A,B,m)$. Then we prove the linear combination property.
\begin{thm}The function $f(z)$ of the form
$(1.1)$ belongs t $M_{w}(A,B,m)$ if and only if
\begin{equation}\sum^{+\infty}_{n=k}n^{m+1}[n(B+1)+A+2]a_{n}\leq 1+A-B .\end{equation}
\end{thm}
\begin{proof} Suppose $(2.1)$ holds and
\begin{equation}
H=|\frac{1}{2}(z-w)(I^{m}f(z))''+(I^{m}f(z))'|-|\frac{B}{2}(z-w)(I^{m}f(z))''+\frac{B+A}{2}|
\end{equation}
Replacing $(I^{m}f(z))'$ and $(I^{m}f(z))''$ by their series
expansion for $0<|z-w|=r<1$, we have
$$H=|\sum^{+\infty}_{n=k}\frac{1}{2}n^{m+1}(n+1)a_{n}(z-w)^{n-1}|-|\frac{1}{2}(B-A-1)\frac{1}{(z-w)^{2}}+\sum^{+\infty}_{n=k}n^{m+1}[\frac{1}{2}(Bn+A+1)]a_{n}(z-w)^{n-1}|$$
$$\leq |\sum^{+\infty}_{n=k}\frac{1}{2}n^{m+1}(n+1)a_{n}r^{n-1}-\frac{1}{2}(1+A-B)\frac{1}{r^{2}}+\sum^{+\infty}_{n=k}n^{m+1}[\frac{1}{2}(Bn+A+1)]a_{n}r^{n-1}|$$
Since this inequality holds for all $r \ \ (0<r<1)$, making
$r\rightarrow 1$, we have
$$H\leq \sum^{+\infty}_{n=k}\frac{1}{2}n^{m+1}[n(B+1)+A+2]a_{n}-\frac{1}{2}(1+A-B)$$
by $(2.1)$
$$ H\leq 0.$$
So we have the required result.\\
Conversely, let\\
$f(z)=\frac{1}{z-w}+\sum^{+\infty}_{n=k}a_{n}(z-w)^{n}$ and
$(1.3)$ holds, so we have
$$|\frac{\sum^{+\infty}_{n=k}\frac{1}{2}n^{m+1}(n+1)a_{n}(z-w)^{n-1}}{\frac{1}{2}(B-A-1)\frac{1}{(z-w)^{2}}+\sum^{+\infty}_{n=k}n^{m+1}[\frac{1}{2}(Bn+A+1)]a_{n}(z-w)^{n-1}}|\leq1,$$
or
$$|\frac{\sum^{+\infty}_{n=k}\frac{1}{2}n^{m+1}(n+1)a_{n}(z-w)^{n+1}}{\frac{1}{2}(1+A+B)-\sum^{+\infty}_{n=k}\frac{1}{2}n^{m+1}(Bn+A+1)a_{n}(z-w)^{n+1}}|\leq 1$$
Since for all $z$,$Re(z)\leq |z|$, so
$$Re\{\frac{\sum^{+\infty}_{n=k}\frac{1}{2}n^{m+1}(n+1)a_{n}(z-w)^{n+1}}{\frac{1}{2}(1+A+B)-\sum^{+\infty}_{n=k}\frac{1}{2}n^{m+1}(Bn+A+1)a_{n}(z-w)^{n+1}}\}\leq 1.$$
Choosing \  $z-w=r$ with $0<r<1$, we get
\begin{equation}\frac{\sum^{+\infty}_{n=k}\frac{1}{2}n^{m+1}(n+1)a_{n}r^{n+1}}{\frac{1}{2}(1+A-B)-\sum^{+\infty}_{n=k}\frac{1}{2}n^{m+1}(Bn+A+1)a_{n}r^{n+1}}\leq
1 .\end{equation} Upon clearing the denominator in $(2.3)$ and
letting $r\rightarrow1$, we get
$$
\sum^{+\infty}_{n=k}\frac{1}{2}n^{m+1}(n+1)a_{n}\leq\frac{1}{2}(1+A-B)-\sum^{+\infty}_{n=k}\frac{1}{2}n^{m+1}(Bn+A+1)a_{n},
$$
or
$$
\sum^{+\infty}_{n=k}\frac{1}{2}n^{m+1}[n(B+1)+A+2]a_{n}\leq\frac{1}{2}(1+A-B).
$$
Hence the proof is complete.
\end{proof}
\begin{thm} Let $f_{j}(z)$ defined by
\begin{equation}
f_{j}(z)=\frac{1}{z-w}+\sum^{+\infty}_{n=k}a_{n,j}(z-w)^{n} ,
j=1,2,...
\end{equation}
be in the class $M_{w}(A,B,m)$, then the function
$$
F(z)=\sum^{t}_{j=1}d_{j}f_{j}(z) \ \ ,\ \ d_{j}\geq0
$$
 is also in $M_{w}(A,B,m)$, where $\sum^{t}_{j=1}d_{j}=1$.
\end{thm}
\begin{proof} Since $f_{j}(z)\in M_{w}(A,B,m)$, by $(2.1)$ we
have,
\begin{equation}
\sum^{+\infty}_{n=k}n^{m+1}[n(B+1)+A+2]a_{n,j}\leq 1+A-B \ \ \
,j=1,2,...
\end{equation}
also
\begin{equation}
F(z)=\sum^{t}_{j=1}d_{j}(\frac{1}{z-w}+\sum^{+\infty}_{n=k}a_{n,j}(z-w)^{n})\nonumber\\
\end{equation}
\begin{equation}
=\frac{1}{z-w}\sum^{t}_{j=1}d_{j}+\sum^{+\infty}_{n=k}(\sum^{t}_{j=1}d_{j}a_{n,j})(z-w)^{n}\nonumber\\
\end{equation}
\begin{equation}
=\frac{1}{z-w}+\sum^{+\infty}_{n=k}s_{n}(z-w)^{n} \ \  where \ \
s_{n}=\sum^{t}_{j=1}d_{j}a_{n,j}\nonumber\\
\end{equation}
But
\begin{equation}
\sum^{+\infty}_{n=k}n^{m+1}[n(B+1)+A+2]s_{n}=\sum^{+\infty}_{n=k}n^{m+1}[n(B+1)+A+2][\sum^{t}_{j=1}d_{j}a_{n,j}]\nonumber\\
\end{equation}
\begin{equation}
=\sum^{t}_{j=1}d_{j}\{\sum^{+\infty}_{n=k}n^{m+1}[n(B+1)+A+2]a_{n,j}\}\nonumber\\
\end{equation}
by $(2.5)$
 $$ \leq\sum^{t}_{j=1}d_{j}(1+A-B)=1+A-B.$$
 Now the proof is complete.
 \end{proof}
\section{Extreme points and Integral operators }
In the last section we investigate about extreme points of
$M_{w}(A,B,m)$ and verify the effect of two operators on functions
in the class $M_{w}(A,B,m)$.
\begin{thm} Let
$$f_{0}(z)=\frac{1}{z-w}   $$and
\begin{equation}f_{n}(z)=\frac{1}{z-w}+\frac{1+A-B}{n^{m+1}[n(B+1)+A+2]}(z-w)^{n}
\ , \ \ n\geq k \end{equation} then $f(z)\in M_{w}(A,B,m) $ if and
only if it can be expressed in the form
$f(z)=\sum_{n=1}^{\infty}C_n f_n(z)$ where $C_{n}\geq0$ ,
$C_{i}=0$ $(i=1,2,...,k-1)$ , $\sum^{+\infty}_{0}C_{n}=1$
\end{thm}
\begin{proof} Let $f(z)=\sum^{+\infty}_{0}C_{n}f_{n}(z)$,
so,
$$f(z)=\frac{C_{0}}{z-w}+\sum^{k-1}_{n=1}C_{n}[\frac{1}{z-w}+\frac{1+A-B}{n^{m+1}[n(B+1)+A+2]}(z-w)^{n}]+$$
$$\sum^{+\infty}_{n=k}C_{n}[\frac{1}{z-w}+\frac{1+A-B}{n^{m+1}[n(B+1)+A+2]}(z-w)^{n}]$$
$$=\frac{1}{z-w}+\sum^{+\infty}_{n=k}\frac{1+A-B}{n^{m+1}[n(B+1)+A+2]}C_{n}(z-w)^{n}.$$
Since
$$\sum^{+\infty}_{n=k}\frac{1+A-B}{n^{m+1}[n(B+1)+A+2]}\frac{n^{m+1}[n(B+1)+A+2]}{1+A-B}C_{n}$$
$$=\sum^{+\infty}_{n=k}C_{n}=\sum^{+\infty}_{n=1}C_{n}=1-C_{0}\leq1
,$$
 so $f(z)\in  M_{w}(A,B,m)$.\\
 Conversely , suppose that $f(z)\in M_{w}(A,B,m)$. Then by  $(2.1)$ we have
 $$0\leq a_{n}\leq \frac{1+A-B}{n^{m+1}[n(B+1)+A+2]}.$$
 By setting
 $$C_{n}=\frac{n^{m+1}[n(B+1)+A+2]}{1+A-B} a_{n}\ \ , \ n\geq 1$$
$C_{i}=0 \ \ (i=1,2,..., k-1) \ ,\
C_{0}=1-\sum^{+\infty}_{n=1}C_{n}$,  we obtain the required
result.
\end{proof}
\begin{thm}Let $\gamma$ be a real number such that $\gamma>1$. If
$f(z)\in M_{w}(A,B,m)$, then the functions
$$H_{1}(z)=\frac{\gamma-1}{(z-w)^{\gamma}}\int^{z}_{w}(t-w)^{\gamma-1}f(t)dt \ \ and$$
$$H_{2}(z)=C\int^{1}_{0}\nu^{C}f(\nu(z-w)+w)d\nu \ , \ C\geq1$$
are also in the same class.\end{thm}
 \begin{proof}Let $f(z)\in M_{w}(A,B,m)$, then a simple
 calculation shows that,
$$H_{1}(z)=\frac{1}{z-w}+\sum^{+\infty}_{n=k}\frac{1}{\gamma+n}a_{n}(z-w)^{n}$$
$$H_{2}(z)=\frac{1}{z-w}+\sum^{+\infty}_{n=k}\frac{C}{C+n+1}a_{n}(z-w)^{n}$$
Since $\frac{C}{C+n+1}<1$ and $\frac{1}{\gamma+n}<1$, by using
Theorem $2.1$, we get,
$$\sum^{+\infty}_{n=k}n^{m+1}[n(B+1)+A+2]a_{n}\frac{1}{\gamma+n}<\sum^{+\infty}_{n=k}n^{m+1}[n(B+1)+A+2]a_{n}<1+A-B$$
and
$$\sum^{+\infty}_{n=k}n^{m+1}[n(B+1)+A+2]a_{n}\frac{C}{C+n+1}<\sum^{+\infty}_{n=k}n^{m+1}[n(B+1)+A+2]a_{n}<1+A-B$$
Hence by Theorem $(2.1)$ we conclude that $H_{1}(z)$ and
$H_{2}(z)$ are in the class $M_{w}(A,B,m)$. So the proof is
complete.
\end{proof}


\bibliography{}

\bibliographystyle{plain}

{\bf Sh. Najafzadeh}\\
Department of Mathematics, Faculty of Science,\\
 University of
Maragheh, Maragheh, Iran \\ E.mail: najafzadeh1234@yahoo.ie \\
{\bf M. Eshaghi Gordji} Department of Mathematics, Semnan
University,\\ P. O. Box 35195-363, Semnan, Iran\\ E.Mail:
madjid.eshaghi@gmail.com\\
{\bf A. Ebadian} \\ Department of
Mathematics, Faculty of Science, \\ Urmia University , Urmia,
Iran \\ E.mail: aebadian@yahoo.com \\
\end{document}